\documentclass[journal,twoside,web]{ieeecolor}
\usepackage{lcsys}

\usepackage{graphicx} 
\usepackage{graphics}
\usepackage{textcomp}
\usepackage{epsfig}
\usepackage{times} 
\usepackage{amsmath}
\usepackage{amssymb}
\usepackage{amsfonts}
\usepackage{float}
\usepackage{siunitx}
\usepackage{scalerel}
\usepackage{cite}
\usepackage{moresize}

\def\BibTeX{{\rm B\kern-.05em{\sc i\kern-.025em b}\kern-.08em
    T\kern-.1667em\lower.7ex\hbox{E}\kern-.125emX}}
\newtheorem{theorem}{Theorem}

\newtheorem{lemma}{Lemma}

 \newtheorem{assumption}{Assumption}
 \newtheorem{Corollary}{Corollary}

\title{\LARGE \bf Distributed Estimation}
\newcommand{\ncom}{\newcommand}

\newcommand{\beqn}{\begin{eqnarray*}}
\newcommand{\eeqn}{\end{eqnarray*}}
\newcommand{\beq}{\begin{eqnarray}}
\newcommand{\eeq}{\end{eqnarray}}

\newcommand{\norm}[1]{\left\lVert #1 \right\rVert}
\newcommand{\inprod}[2]{\left\langle #1, #2 \right\rangle}
\ncom\R{\mathbb{R}}
\DeclareMathOperator*{\argmin}{arg\,min}

\author{Anik Kumar Paul, Arun D Mahindrakar and Rachel K Kalaimani
\thanks{Anik  is a Graduate student in the Department of Electrical Engineering, IIT Madras, Chennai-600036, India
        {email: anikpaul42@gmail.com}}
        \thanks{Arun  and Rachel are  with the Department of Electrical Engineering, Indian Institute of Technology Madras, Chennai-600036, India (email: arun\_dm@iitm.ac.in, rachel@ee.iitm.ac.in) }
    }

\begin{document}
\title{Almost Sure Convergence and Non-asymptotic Concentration Bounds for Stochastic Mirror Descent Algorithm}
\maketitle
\thispagestyle{plain}
\pagestyle{plain}
\begin{abstract}
This letter investigates the convergence and concentration properties of the Stochastic Mirror Descent (SMD) algorithm utilizing biased stochastic subgradients. We establish the almost sure convergence of the algorithm's iterates under the assumption of diminishing bias. Furthermore, we derive concentration bounds for the discrepancy between the iterates' function values and the optimal value, based on standard assumptions. Subsequently, leveraging the assumption of Sub-Gaussian noise in stochastic subgradients, we present refined concentration bounds for this discrepancy.
\end{abstract}
\vspace{-1em}
\section{Introduction}
The mirror descent algorithm extends gradient descent to non-Euclidean spaces, facilitating optimization in complex geometries by using proximal functions like the Bregman divergence \cite{paul2024convergenceanalysisstochasticsaddle,anik,7004065}. Recently, it has gained attention in large-scale optimization, data-driven control, power systems, robotics, and game theory \cite{8409957}. This letter focuses on the SMD algorithm, crucial for scenarios involving robust estimation, resource allocation, and learning theory \cite{1055386,KellyMT98,811451}. The SMD algorithm overcomes the impracticality of exact subgradients by employing noisy subgradients, which are particularly effective in these scenarios.
 Extending the rich history of noisy gradients in gradient descent algorithms (see \cite{Polyak87}), the SMD algorithm was developed in \cite{doi:10.1137/120894464}. There has been a surge of interest in various stochastic subgradient and mirror descent algorithm variants \cite{JMLR:v23:20-899,7405263,articlesmd}, usually assuming zero-mean noise, leading to unbiased stochastic subgradients.  However, this assumption often fails in practical applications, especially in zeroth-order optimization, where stochastic noise typically has a non-zero mean, resulting in biased stochastic subgradients \cite{doi:10.1137/120880811,Nesterov2015RandomGM}.

In this letter, we analyze the SMD algorithm with non-zero bias. The study of stochastic gradient descent and its variants with non-zero bias has been actively researched due to its broad applications in machine learning, reinforcement learning, and signal processing \cite{9186148}. Recent works, such as \cite{ajalloeian2021convergence} for gradient descent and \cite{9416872} for mirror descent, mainly focus on bounding the expected error between the function value and the optimal value. While expected error rates are important, it is also crucial to determine if the iterates of a single trial of the algorithm converge to the optimal solution. The work in \cite{10143924} advances this analysis by demonstrating the almost sure convergence of iterates for the mirror descent algorithm with non-zero bias.

The studies in \cite{10143924} showed that iterates converge to a neighborhood around the optimal value due to non-zero but bounded bias. In this letter, we demonstrate that under certain conditions on the bias, the iterates of the SMD algorithm almost surely converge to the optimal solution. We also determine the almost sure convergence rate, estimating the number of iterations needed to ensure the function value of the iterates is almost surely below a certain threshold.
\\ Additionally, we derive concentration inequalities for the SMD algorithm with non-zero diminishing bias, providing probability bounds for the error under assumptions of bounded variance and Sub-Gaussian distribution of the stochastic subgradient.  This area is less explored compared to expected error bounds, with recent progress focusing on stochastic gradient descent with unbiased gradients (see \cite{JMLR:v23:21-0560} and references therein). We develop the concentration inequality  for  the SMD algorithm, particularly focusing on biased subgradients. We provide the bounds on the  probability of the error between function value  exceeding a specified threshold under two assumptions: bounded variance of the stochastic gradient and Sub-Gaussian distribution of the stochastic subgradient. Our analysis indicates the minimum number of iterations required to ensure the error remains below a given threshold with a specified confidence level. Importantly, to the best of our knowledge, such concentration inequalities for biased stochastic subgradients have not been previously explored.

This novel analysis reveals the SMD algorithm's performance over finite periods, emphasizing how bias impacts its efficacy and showing how selecting an appropriate Bregman divergence can enhance performance compared to the standard subgradient method.

The contributions of this letter are listed as follows:
\begin{enumerate}
    \item We analyze the SMD algorithm with non-zero diminishing bias. We demonstrate that, under the summability condition on product of step=size and bound on bias, the iterates almost surely converge to the set of optimal solutions. Additionally, we derive the almost surely convergence rate, achieving the optimal rate of $\mathcal{O}(\frac{1}{\sqrt{t}})$ with suitable step-sizes.
    \item  Furthermore, we derive concentration inequalities for the SMD algorithm with non-zero diminishing bias. We provide bounds under two scenarios: standard assumptions on stochastic subgradients and the more restrictive Sub-Gaussian distribution assumption. This result aids in calculating the number of iterations required for any desired confidence level and specified neighborhood, ensuring that the function values of the iterates lie within the specified neighborhood of the optimal function value with that confidence after those iterations. We also outline sufficient conditions that validate the Sub-Gaussian assumption, highlighting its impact on improving concentration bounds. 
\end{enumerate}
\vspace{-1.2em}
\subsection{Notation and Math Preliminary}
\label{mathpre}
Let $\R$ and $\mathbb{R}^n$  represent the set of real numbers, set of $n$ dimensional real vectors.   Let~$\vert\vert \cdot
\vert\vert$ denote any \mbox{norm} on~$\R^n$.  Given a norm $\norm{.}$ on $\mathbb{R}^n$, the dual norm of $x\in \R$ is $\norm{x}_\ast := \sup\{ \inprod{x}{y}: {\norm{y}\leq 1}, y\in \R^n \}$, where $\inprod{x}{y}$ denotes the standard inner-product on $\R^n$. For $x\in \mathbb{R}^n$, $[x]_i$ denotes the $i$th component of the vector $x$. The set of all natural numbers is denoted by $\mathbb{N}$ and $I_N\in \R^{N\times N}$ is the identity matrix. The set of integers from $1$ to $N$ is denoted by $[N]$. 
A random vector  $X\sim \mathcal{N}(0_n, I_n)$ denotes a $n$-dimensional normal random vector with zero-mean and unit standard-deviation. Let $(\{X_i\}, \; 1 \leq i \leq n )$ be the sequence of random vectors defined on the probability space $(\Omega,\mathcal{F},\mathbb{P})$ then, $\sigma(\{X_i\}, \; 1 \leq i \leq n )$ denotes the smallest sigma algebra generated by the random vectors $(\{X_i\}, \; 1 \leq i \leq n)$.

Let  $f: \mathbb{R}^n \to \mathbb{R}$ be a  convex function.  A vector $g\in \R^n$  is called a subgradient of $f$ at $x$ if and only if $f(y)\geq f(x)+\inprod{g}{y-x}\; \forall \; y\in \R^n $. 
 An event $A \in \mathcal{F}$ is occurred almost surely (a.s.) if $\mathbb{P}(A) =1$. A random variable $X \in \mathcal{L}^1$ if $\mathbb{E}[\norm{X}] < \infty$. A sequence of random variables $\{X_i\}_{i\geq 1} \subseteq \mathcal{L}^1$ converges to some random variable $X^\ast$ almost surely if $\mathbb{P}(\{\omega \in \Omega|\lim\limits_{i \to \infty}X_i(\omega) \to X^\ast(\omega)\}) = 1$. For the sequence $\{X_i\}_{i \geq 1}$, we say $X_i = \mathcal{O}(a_i)$ a.s. $(\{a_i\} \subseteq \mathbb{R})$ if $\limsup\limits_{i \to \infty} \frac{|X_i|}{a_i} < \infty$ a.s, in other words, $\exists \; C_1 \in \mathcal{L}^1$ such that $|X_i|\leq C_1 a_i$ a.s \cite{https://doi.org/10.48550/arxiv.1108.3924}.
 The following theorem will be helpful in our analysis .
\begin{theorem}{\cite{ROBBINS1971233}}
Let ${X(t)}$, ${Y(t)}$, and ${Z(t)}$ be non-negative random processes adapted to a filtration ${\mathcal{F}_t}$, such that   $\mathbb{E}[X(t+1)| \mathcal{F}t] \leq (1 + \alpha(t)) X(t) - Y(t) + Z(t)$ a.s.
where $\alpha(t) \in \mathbb{R}$ satisfies $\sum\limits_{t \geq 1} \alpha(t) < \infty$. Additionally, suppose  $\sum\limits_{t \geq 1} Z(t) < \infty$ almost surely. Then  $X(t)$ converges to a non-negative $\mathcal{L}^1$ random variable  and $\sum\limits_{t \geq 1} Y(t) < \infty$ a.s. 
\label{Robbins-Siegmund}
\end{theorem}
\vspace{-1.43em}
\section{Problem Statement}
In this letter, we consider the following optimization problem 
\vspace{-1em}
\begin{equation*}\tag{CP1}
    \min\limits_{x \in \mathbb{X}} f(x).
    \label{CP1}
\end{equation*}
The constraint set $\mathbb{X} \subseteq \mathbb{R}^n$ is  convex and compact with diameter $D$.  The function $f:  \mathbb{R}^n \to \mathbb{R}$ is a convex and continuous.   Define $ f^\ast = \min\limits_{x \in \mathbb{X}} f(x) \; \; \text{and} \; \;  \mathbb{X}^\ast = \{x^\ast \in \mathbb{X} | f(x^\ast) = f^\ast\}.$
The set  $\mathbb{X}^\ast$ is nonempty due to compactness of the constraint set $\mathbb{X}$ and continuity of the function $f$.
 Since  $\mathbb{X}$ is compact, subgradient of  $f$ is bounded. In other words, there exists a constant  $G > 0$ such that $\forall \; g \in \partial f(x)$ and $ \; \forall x \in \mathbb{X}$, we have $\norm{g}_\ast \leq G$.
\vspace{-0.9em}
\subsection{SMD Algorithm}
Before moving to  the algorithmic steps, we  define
the notion of Bergman divergence. Let $R$ be the $\sigma_R$-strongly convex function and differentiable over an open set that contains the set $\mathbb{X}$. The  Bergman divergence $\mathbb{D}_R(x,y) : \mathbb{X} \times \mathbb{X} \to \mathbb{R}$ is 
$ \mathbb{D}_{R}(x,y) := R(x)- R(y) - \inprod{\nabla R(y)}{x-y} \; \; \forall \; \; x,y \in \mathbb{X}.$
It is clear from the definition of strong convexity that
\\ $\mathbb{D}_R(x,y) \geq \frac{\sigma_R}{2}\norm{x-y}^2$ and $\forall \; \;  x,y,z \in  \mathbb{X}$
\\ $\mathbb{D}_R(z,y)- \mathbb{D}_R(z,x)-\mathbb{D}_R(x,y) = \inprod{\nabla R(x) - \nabla R(y)}{z-x}.$
At iteration $t$, let $x(t)$ be the iterates of the SMD algorithm. 
The next iterates at $(t+1)$ is calculated as follows:
\begin{equation}
    x(t+1) = \argmin\limits_{x \in \mathbb{X}} \{ \inprod{\Tilde{g}(t)}{x-x(t)} + \frac{1}{\alpha(t)} \mathbb{D}_R(x,x(t))\}
    \label{jomli}
\end{equation}
where, $\Tilde{g}(t)$ denotes the noisy subgradient of the function $f(x)$ at $x = x(t)$. The parameter $\alpha(t)$  represents the step-size of the algorithm. 
The ergodic averages of  $x(t)$ is as follows: 
 $    z(t) = A(t)^{-1}\sum\limits_{j=1}^{t} \alpha(j)x(j)$, where, $A(t) = \sum\limits_{k=1}^{t} \alpha(k)$.
In this context, let us define the filtration $\{ \mathcal{F} \}_{t \geq 1}$ as follows:
$    \mathcal{F}_t = \sigma ( \{x_l\} | 1 \leq l \leq t )  \; \; \forall \; t \in \mathbb{N} $.

Define another filtration as $\{\mathcal{G}_t\}_{t \geq 1}$ such that $\mathcal{G}_{t-1} = \mathcal{F}_t$, which will be helpful in the subsequent analysis. Before proceeding further, we need the following Assumptions on stochastic subgradient and the step-size.  
\begin{assumption}
 $\mathbb{E} [\Tilde{g}(t) | \mathcal{G}_{t-1}] = g(t) + b(t)$, $ \mathbb{E}[\norm{\Tilde{g}(t)}_\ast^2 | \mathcal{F}_t] \leq \nu^2 \; \text{a.s.}$
    \label{ass1}
\end{assumption}
\begin{assumption} 
    $ \sum\limits_{t \geq 1}  \alpha(t) = \infty \; \; \text{and} \; \; \sum\limits_{t \geq 1} \alpha(t)^2 < \infty.$
    \label{ass2}
\end{assumption}
\begin{assumption} Let $B(t) = \norm{b(t)}_\ast$ ,   $\sum\limits_{t \geq 1} \alpha(t) B(t) < \infty.$
    \label{ass3}
\end{assumption}
Assumption \ref{ass3} regarding bias may seem restrictive, but it is often valid in many zeroth-order optimization applications. For instance, if $\Tilde{g}(t)$ is computed using the Gaussian approximation as described in \cite{Nesterov2015RandomGM}, choosing the smoothing parameter $\mu(t) = \mathcal{O}(\frac{1}{t})$ ensures that Assumption \ref{ass3} is satisfied.
 One important consequence of Assumptions \ref{ass1} and \ref{ass3} is that the stochastic subgradient can be rewritten as follows:
 $  \Tilde{g}(t) = g(t) + b(t) + \zeta(t)$.
where, $b(t)$ is a $\mathcal{F}_t$ measurable random variable and $\mathbb{E}[\zeta(t) | \mathcal{F}_t] = 0$ a.s.
Under these Assumptions, the ensuing section aims to address the following questions: $1)$  Whether the sequence of random vectors $\{x(t)\}_{t \geq 1}$ converges  to $\mathbb{X}^\ast$ a.s., and   if so  what is the convergence rate? $2)$ What should be the concentration bound for  $f(z(t)) - f^\ast$ ?
\vspace{-0.6em}
\section{Performance Analysis of SMD Algorithm}
In this section, we focus on demonstrating the almost sure convergence of the iterates  for both $x(t)$ and $z(t)$. 
\vspace{-1.1em}
\subsection{Almost Sure Convergence}
\begin{theorem}
    The sequence of iterates $\{x(t)\}$ generated by the SMD algorithm, as described in  \eqref{jomli}, almost surely converges to   $x^\ast \in \mathbb{X}^\ast$.
    \label{thehe}
\end{theorem}
\begin{proof}
    The  first-order optimality condition to  \eqref{jomli} yields
\begin{equation}
\begin{split}
     & \alpha(t) \inprod{\Tilde{g}(t)}{x - x(t+1)} \; \; \; (\forall \; x \in \mathbb{X})
        \\  \geq &   -  \inprod{\nabla R (x(t+1))- \nabla R (x(t))}{x-x(t+1)}
        \\ = & \mathbb{D}_R(x{(t+1)},x(t))+ \mathbb{D}_R(x,x{(t+1)})-\mathbb{D}_R(x,x(t)).
          \end{split}
          \label{firstor}
\end{equation}
 From the LHS of \eqref{firstor}, we obtain
\begin{equation*}
    \begin{split}
        &\alpha(t) \inprod{\Tilde{g}(t)}{x - x{(t+1)}}  \leq \alpha(t) \inprod{\Tilde{g}(t)}{x-x(t)} +
        \\ &    \frac{\alpha(t)^2}{2 \sigma_R} \norm{\Tilde{g}(t)}_\ast^2 + \frac{\sigma_R}{2} \norm{x(t)-x{(t+1)}}^2.
    \end{split}
\end{equation*}
The last inequality  follows by applying the Young-Fenchel inequality to the term $\alpha(t)  \inprod{\Tilde{g}(t)}{x(t) - x(t+1)}$. Hence from \eqref{firstor}, we get that $ \mathbb{D}_R(x,x{(t+1)})$
\begin{equation}
    \begin{split}
         \leq & \mathbb{D}_R(x,x(t)) + \alpha(t) \inprod{\Tilde{g}(t)}{x-x(t)} +  \frac{\alpha(t)^2}{2 \sigma_R} \norm{\Tilde{g}(t)}_\ast^2.
    \end{split}
    \label{Ber}
\end{equation}
Notice that $\mathbb{D}_R(x{(t+1)},x(t)) \geq \frac{\sigma_R}{2} \norm{x{(t+1)}-x(t)}^2$.

Letting $x = x^\ast \in \mathbb{X}^\ast$ in \eqref{Ber} and considering only the term $\alpha(t) \inprod{\Tilde{g}(t)}{x^\ast-x(t)}$ and upon taking conditional expectation on both sides, we get $\mathbb{E}[\mathbb{D}_R
(x^\ast,x(t+1))| \mathcal{F}_t]$
\begin{flalign}
          \leq  \mathbb{D}_R(x^\ast,x(t)) + \alpha(t) \inprod{g(t) + b(t)}{x^\ast -x(t)} + \frac{\alpha(t)^2}{2 \sigma_R} \nu^2. && \nonumber
    \label{gd}
\end{flalign}
This inequality holds true in view of Assumption \ref{ass1}. Consider the term $\alpha(t) \inprod{g(t) + b(t)}{x^\ast -x(t)} $ 
\begin{equation*}
    \begin{split}
         \leq & \alpha(t) (f(x^\ast)-f(x(t)) + \alpha(t)   (\norm{x^\ast-x(t)}^2 + 1) \norm{b(t)}_\ast
    \end{split}
    \label{rey}
\end{equation*}
\begin{equation*}
    \leq  \alpha(t) (f(x^\ast)-f(x(t)) + \alpha(t) B(t) + \frac{2 \alpha(t) B(t)}{\sigma_R}( \mathbb{D}_R(x^\ast,x(t))) .
\end{equation*}
The first inequality in the above equation arises from the convexity of the function $f$.  By applying the generalized Cauchy-Schwarz inequality, we obtain 
\\ $\inprod{b(t)}{x^\ast-x(t)} \leq \norm{b(t)}_\ast \norm{x^\ast-x(t)} \leq (\norm{x^\ast-x(t)}^2 +1 ) B(t) $
 The subsequent inequality  follows from strong convexity of $R$.
Hence, we get that 
\begin{flalign}
        & \mathbb{E}[\mathbb{D}_R(x^\ast,x(t+1))| \mathcal{F}_t] \leq (1+ \frac{2 \alpha(t) B(t)}{\sigma_R})\mathbb{D}_R(x^\ast,x(t)) \nonumber
        \\ & - \alpha(t) (f(x(t))-f^\ast) + \alpha(t)B(t) + \frac{\alpha(t)^2}{2 \sigma_R} \nu^2.
    \label{rs}
\end{flalign}
\\ Notice that $\sum\limits_{t \geq 1} \alpha(t) B(t) < \infty$ and $\sum\limits_{t \geq 1} \alpha(t)^2 < \infty$ in view of Assumption \ref{ass2} and \ref{ass3}. Hence, by applying Theorem \ref{Robbins-Siegmund} it follows that  $\mathbb{D}_R(x^\ast,x(t))$ converges to a $\mathcal{L}^1$ random variable $\forall \; x^\ast \in \mathbb{X}^\ast$ a.s. and   $\sum\limits_{t \geq 1} \alpha(t) (f(x(t))-f^\ast) < \infty$ a.s.
\\ Notice that  it implies that $  \liminf\limits_{t \to \infty} (f(x(t))-f^\ast) = 0$.
Otherwise, $\forall \; \epsilon > 0$ $\exists \; t_0 > 0$ such that $\forall \; t \geq t_0$ we have $f(x(t))-f^\ast \geq \epsilon$. That implies $ \sum\limits_{t \geq t_0} \alpha(t)(f(x(t))-f^\ast) \geq \epsilon \sum\limits_{t \geq t_0} \alpha(t) = \infty,$ a contradiction.
Hence, there exists a subsequence $\{x(t_k)\}$ of $\{x(t) \}$ such that $\lim\limits_{k \to \infty} f(x(t_k)) = f^\ast$. 

However, since $x(t_k) \in \mathbb{X}$, a compact set,  there exists a convergent subsequence. Without loss of generality, we assume that $x(t_k) \to x_1 \in \mathbb{X}$. Due to the continuity of $f$, we have $f(x(t_k)) \to f(x_1)$. Thus, we conclude that $x_1 \to \mathbb{X}^\ast$. 

This implies that $\mathbb{D}_R(x_1,x(t_k))$ converges to zero almost surely. Since it has already been proven that  $\mathbb{D}_R(x^\ast,x(t))$ converges a.s,  we can conclude that $x(t) \to x_1 \in \mathbb{X}^\ast$.
\end{proof}
The next Corollary is true in view of Toeplitz Theorem  \cite[Theorem 5, pp. 75]{Knopp1990TheoryAA}.
\begin{Corollary}
   The sequence of iterates $\{z(t)\}$  also almost surely converges to the set $\mathbb{X}^\ast$.
\end{Corollary}
In the remaining part of this subsection and the next subsection, we focus on the sequence $\{z(t)\}$ to provide the almost sure convergence rate and to conduct the finite-time analysis.
\begin{theorem}
    The almost sure convergence rate of the SMD algorithm, in terms of the function value, is given by  $f(z(t)) - f^\ast = \mathcal{O} \Big{(}\frac{1}{\sum\limits_{j=1}^{t} \alpha(j)}\Big{)}$ a.s.
    \label{almstrate}
\end{theorem}
\begin{proof}
    Note that $z(t)$ can be written also in the following way
     $ z(t) = \beta(t) x(t) + (1 - \beta(t)) z(t-1)$
    where, $\beta(t)$ is computed as follows $\beta(t) = A(t)^{-1} \alpha(t) \in (0,1)$.
From Jensen's inequality
 $ f(z(t)) \leq \beta(t) f(x(t)) + (1- \beta(t)) f(z(t-1))$.
    Hence, we have  $\alpha(t)f(x(t))
        \geq \sum\limits_{j=1}^{t} \alpha(j) f(z(t)) - \sum\limits_{j=1}^{t-1} \alpha(j) f(z(t-1)) $.
      From \eqref{rs} of Theorem \ref{thehe}, we obtain $\mathbb{E}[\mathbb{D}_R(x^\ast,x(t+1)) + \sum\limits_{j=1}^{t} \alpha(j)(f(z(t))-f^\ast)| \mathcal{F}_t] \leq $
     \vspace{-0.7em}
    \begin{equation*}
    \begin{split}
          &  (1+ \frac{2 \alpha(t) B(t)}{\sigma_R})(\mathbb{D}_R(x^\ast,x(t)) + \sum\limits_{j=1}^{t-1} \alpha(j) (f(z(t-1))-f^\ast)
        \\ &  + \alpha(t)B(t) + \frac{\alpha(t)^2}{2 \sigma_R} \nu^2.
    \end{split}
\end{equation*}
However, by applying  Theorem \ref{Robbins-Siegmund} we get that
$\mathbb{D}_R(x^\ast,x(t)) + \sum\limits_{j=1}^{t-1} \alpha(j) (f(z(t))-f^\ast)$ converges to some $\mathcal{L}^1$ random variable. 
It is already shown in Theorem \ref{thehe} that $\mathbb{D}_R(x^\ast,x(t))$ converges to a $\mathcal{L}^1$ random variable. \\ This implies w.p. $1$  $\limsup\limits_{t \to \infty} \sum\limits_{j=1}^{t} \alpha(j) (f(z(t)-f^\ast) < \infty.$
\end{proof}
   Theorem \ref{almstrate} reveals the almost sure convergence rate of the SMD algorithm based on step-size.  For $\alpha(t) = \mathcal{O}(\frac{1}{t^k})$ with $k \in (\frac{1}{2}, 1)$, the algorithm achieves an almost sure convergence rate of $\mathcal{O}(\frac{1}{t^{1-k}})$, approaching the optimal rate of $\mathcal{O}(\frac{1}{\sqrt{t}})$. These results align with the standard almost sure convergence rate for stochastic gradient descent with unbiased gradients \cite{pmlr-v134-sebbouh21a}. However, Theorem \ref{almstrate} does not address the algorithm's performance at specific finite times or the dependence of the convergence rate on the bias term. In the next subsection, we analyze finite-time performance and explore the convergence rate in terms of various  other parameters. 
   \vspace{-1.2em}
   \subsection{Non asymptotic Concentration Bound}
   \begin{theorem}
       For any $\epsilon > 0$, $\exists \; t_0 \in \mathbb{N}$ such that
        $\forall \; t \geq t_0$ we have $\mathbb{P}(f(z(t))-f^\ast \geq \epsilon )$
       \begin{equation*}
           \begin{split}
             \leq & \frac{3\nu^2\sum\limits_{k=1}^{t}\alpha(k)^2}{2\sigma_R\epsilon \mathrm{K} \sum\limits_{k=1}^{t} \alpha(k)}  + \frac{9 D^2 \kappa_1 \sum\limits_{k=1}^{t} \alpha(k)^2(\nu^2 + G^2 + B(k)^2)}{\epsilon^2 \mathrm{K}^2 \Big{(}\sum\limits_{k=1}^{t} \alpha(k) \Big{)}^2}.
           \end{split}
       \end{equation*}
       \label{finite}
   \end{theorem} 
    \begin{proof}
      Based on the proof of Theorem \ref{thehe}, from  \eqref{rs} , we obtain the following inequality:
       \begin{flalign}
        & \mathbb{D}_R(x^\ast,x(t+1)) \leq (1+ \rho(t))\mathbb{D}_R(x^\ast,x(t)) + \frac{\alpha(t)^2}{2 \sigma_R} \norm{\Tilde{g}(t)}_\ast^2  && \nonumber
        \\ & - \alpha(t) ( (f(x(t))-f^\ast) - \inprod{\zeta(t)}{x^\ast-x(t)} - B(t)) &&
    \label{try}
\end{flalign}
where, $\rho(k) = \sigma_R^{-1} 2 \alpha(k) B(k)$
and define three new random variables $Z(t)$, $X(t)$ $W(t)$ as
\\ $Z(t) = \mathbb{D}_R(x^\ast,x(t)) \times \prod\limits_{k=1}^{t-1} (1 + \rho(k) )^{-1}$, 
\\ $ W(t) = \alpha(t)\inprod{\zeta(t)}{x^\ast-x(t)} \times \prod\limits_{k=1}^{t} (1 + \rho(k))^{-1}$ and 
\\ $ X(t) = \alpha(t)(f(x(t))-f^\ast)\times  \prod\limits_{k=1}^{t} (1 + \rho(k))^{-1}$. 
From the definition of  $X(t)$, we observe that \\ $X(t)$
 $ \geq \alpha(t)(f(x(t))-f^\ast)  \prod\limits_{k=1}^{\infty} (1 + \rho(k))^{-1}$
 \\  \hspace{2em}  $ \geq \alpha(t)(f(x(t))-f^\ast)   \prod\limits_{k=1}^{\infty} \exp{(- \rho(k))}$
\vspace{-1.2em}
\begin{flalign}
      = \alpha(t)(f(x(t))-f^\ast)   \exp{(-\sum\limits_{k=1}^{\infty} \rho(k))  }. &&
\label{func}
   \end{flalign}
  The second inequality in \eqref{func} holds because of the inequality $1+x \leq e^x$ $\forall \; x > 0$.  Define $\mathrm{K} =  \exp{(-\sum\limits_{k=1}^{\infty} \rho(k))  } $.
  Therefore, multiplying both sides of \eqref{try} by $\prod\limits_{k=1}^{t} (1+ \rho(k))^{-1}$, we get 
  \\ $ Z(t+1)  \leq Z(t) + \alpha(t)B(t) + \frac{\alpha(t)^2}{2 \sigma_R} \norm{\Tilde{g}(t)}_\ast^2
            -\alpha(t) \mathrm{K} (f(x(t))-f^\ast) + W(t).$
  Applying the telescopic sum and rearranging the equation, we get 
 $\mathrm{K} \sum\limits_{k=1}^{t} \alpha(k) (f(x(k))-f^\ast)  \leq Z(1) $
  \vspace{-0.8em}
 \begin{flalign}
             + \sum\limits_{k=1}^{t} W(k)  + \sum\limits_{k=1}^{t} \Big{(}\frac{\alpha(k)^2}{2 \sigma_R} \norm{\Tilde{g}(k)}_\ast^2 + \alpha(k) B(k) \Big{)}. 
     \label{gdggd}
 \end{flalign}
 From Jensen's inequality  
\\  $ \sum\limits_{k=1}^{t} \alpha(k) f(z(t)) \leq \sum\limits_{k=1}^{t} \alpha(k) f(x(k)).$
 Hence, from \eqref{gdggd}, we get that $  f(z(t))-f^\ast   \leq \frac{1}{\mathrm{K} \sum\limits_{k=1}^{t} \alpha(k)} \times $
 \vspace{-1em}
 \begin{flalign}
     \Big{(} Z(1)  + \sum\limits_{k=1}^{t} \alpha(k) B(k) \nonumber   + \sum\limits_{k=1}^{t} W(k) + \sum\limits_{k=1}^t \frac{\alpha(k)^2}{2 \sigma_R} \norm{\Tilde{g}(k)}_\ast^2 \Big{)}. &&
 \end{flalign}
  Note that $Z(1) = \mathbb{D}_R(x^\ast,x(1))$. 
 Consider $t_0 \in \mathbb{N}$ such that $\forall \; t \geq t_0$, we have $ \mathbb{D}_R(x^\ast,x(1)) + \sum\limits_{k=1}^{t} \alpha(k) B(k)  \leq  \tau(t) $, where, $\tau(t) = \frac{\epsilon}{3} \mathrm{K} \sum\limits_{k=1}^{t} \alpha(k)$.
Such a $t_0$ always exists and is finite in view of Assumption \ref{ass2} and \ref{ass3}. 
 Consider any $t \geq t_0$ and suppose $f(z(t))-f^\ast \geq \epsilon$. Then   either  $\sum\limits_{k=1}^{t} W(k) \geq \tau(t)$ or, 
  $\sum\limits_{k=1}^{t} \frac{\alpha(k)^2}{2 \sigma_R} \norm{\Tilde{g}(k)}_\ast^2 \geq \tau(t) $ holds. That is, 
  \begin{flalign}
      \mathbb{P}\Big{(}f(z(t))-f^\ast \geq \epsilon \Big{)} &&
      \label{24}
  \end{flalign}
  \vspace{-2em}
\begin{flalign}
   \leq  \mathbb{P} \Big{(}\sum\limits_{k=1}^{t} W(k) \geq \tau(t) \Big{)}   + \mathbb{P}\Big{(}\sum\limits_{k=1}^{t} \frac{\alpha(k)^2}{2 \sigma_R} \norm{\Tilde{g}(k)}_\ast^2 \geq \tau(t) \Big{)}. && \nonumber
   \end{flalign}
We apply Markov's inequality to get a bound on the RHS. 
Hence, $\mathbb{P}\Big{(}\sum\limits_{k=1}^{t} \frac{\alpha(k)^2}{2 \sigma_R} \norm{\Tilde{g}(k)}_\ast^2 \geq \tau(t) \Big{)}$ 
\vspace{-0.5em}
\begin{equation}
    \begin{split}
  \leq & \frac{\mathbb{E}\Big{[} \sum\limits_{k=1}^{t} \frac{\alpha(k)^2}{2 \sigma_R} \norm{\Tilde{g}(k)}_\ast^2\Big{]}}{\tau(t)}
 \leq  \frac{\nu^2\sum\limits_{k=1}^{t}\alpha(k)^2}{2\sigma_R\epsilon \tau(t) }.
    \end{split}
    \label{mama}
\end{equation}
For the other term we get
\\ $ \mathbb{P} \Big{(}\sum\limits_{k=1}^{t} W(k) \geq \tau(t) \Big{)}  =    \mathbb{P} \Big{(}(\sum\limits_{k=1}^{t} W(k) )^2 \geq \tau(t)^2 \Big{)}$.
\vspace{-0.4em}
Consider $W(t) = \alpha(t)  \prod\limits_{k=1}^{t} (1 + \rho(k))^{-1} \inprod{\zeta(t)}{x^\ast-x(t)} $, 
\\ Therefore,
$ W(t)^2 \leq \alpha(t)^2 \norm{\zeta(t)}_\ast^2 \norm{x^\ast-x(t)}^2$. 
\\ Note that $\prod\limits_{k=1}^{t} (1 + \rho(k))^{-1} < 1$
By taking the expectation on both sides and considering that the diameter of the constraint set is $D$, we get 
\\ $ \mathbb{E}[ W(t)^2] \leq \alpha(t)^2 D^2\mathbb{E}[\norm{\zeta(t)}_\ast^2] \leq \alpha(t)^2 D^2 \kappa_1 (\nu^2 + G^2 + B(t)^2) .$ The last inequality follows from the fact $ \mathbb{E}[\norm{\zeta(t)}_\ast^2]  = \mathbb{E}[\norm{\Tilde{g}(t)-g(t)-b(t)}_\ast^2]$, where we have used parallelogram law and $\kappa_1 >0$ is a norm equivalence constant.

Now consider any $i, j \in \mathbb{N}$ such that $i < j$ then 
\\ $\mathbb{E}[W(i)W(j)] = \mathbb{E}[W(i)\mathbb{E}[W(j) | \mathcal{F}_j]] = 0.$
This equality holds from the fact $W(i)$ is $\mathcal{F}_j$ measurable and 

$\mathbb{E}[W(j)| \mathcal{F}_j] = \alpha(j) e(j) \mathbb{E}[\inprod{\zeta(j)}{x^\ast-x(j)}| \mathcal{F}_j] = 0$.
 Hence, from Markov's inequality $ \mathbb{P} \Big{(}\sum\limits_{k=1}^{t} W(k) \geq \tau(t) \Big{)} \leq  $
  \vspace{-0.68em}
 \begin{flalign}
      \frac{ D^2 \kappa_1 \sum\limits_{k=1}^{t} \alpha(k)^2(\nu^2 + G^2 + B(k)^2)}{\tau(t)^2}. &&
      \label{modm}
  \end{flalign}
   Adding the above equation with \eqref{mama}, we get the result. 
   \end{proof}
  Theorem \ref{finite} quantifies the probability of the SMD algorithm's function value deviating from the optimal value at a finite time and estimates the convergence speed based on factors like the Lipschitz constant, strong convexity parameter of $R$, and choice of $B(t)$. The following corollary further clarifies these aspects.  
\begin{Corollary}
For any $\epsilon > 0$ and a confidence level $0 < p < 1$, let $p_1 = 1-p$.   Let $\mathcal{R} = \sup\limits_{(x,y) \in \mathbb{X}} \mathbb{D}_R(x,y)$.
Define $t_0$  such that $\forall \;  t \geq t_0$,
\\ $ \sum\limits_{k=1}^{t} \alpha(k) \geq \frac{3}{\epsilon} \exp{\Big{(}\sum\limits_{k=1}^{\infty} \frac{2 \alpha(k)B(k)}{\sigma_R}\Big{)}} \Big{(} \mathcal{R} + \sum\limits_{k=1}^{t} \alpha(k) B(k)\Big{)}$
\\  Define $t_1$ such that $\forall \; t \geq t_1$,
\\ $ \sum\limits_{k=1}^{t} \alpha(k) \geq \frac{3 \nu^2}{\sigma_R \epsilon p_1} \exp{\Big{(}\sum\limits_{k=1}^{\infty} \frac{2 \alpha(k)B(k)}{\sigma_R}\Big{)}} \sum\limits_{k=1}^{t} \alpha(k)^2$
 \\ and $t_2$ such that $\forall \; t \geq t_2$,
\\ $ (\sum\limits_{k=1}^{t} \alpha(k))^2 \geq \frac{9D^2 \kappa_1}{\epsilon^2}  \exp{\Big{(}\sum\limits_{k=1}^{\infty} \frac{2 \alpha(k)B(k)}{\sigma_R}\Big{)}} \times
    \sum\limits_{k=1}^{t} \alpha(k)^2 (\nu^2 + G^2 + B(k)^2)$.
Then $\forall \; t$ $ \geq$ $ \max \{ t_0, t_1, t_2\}$, we have
\\ $\mathbb{P}( f(z(t)) -f^\ast <  \epsilon ) \geq p.$
Notice that  $\max \{ t_0, t_1, t_2\}$  $< \infty$ due to Assumptions \ref{ass2} and \ref{ass3}.
\label{concentrationc}
\end{Corollary}
In this subsection, we detail how the iterations of the SMD algorithm converge to the optimal solution, dependent on various parameters. The concentration bound in Theorem \ref{thehe} can be significantly enhanced under a stricter assumption: that the stochastic subgradient is Sub-Gaussian, which is the focus of the next subsection. 
\vspace{-1em}
\subsection{Non Asymptotic Concentration Bound under Sub-Gaussian Assumptions}
We begin this subsection with the assumption that the noise in the stochastic gradient, $\zeta(t)$, is a Sub-Gaussian random vector. This assumption is clarified in the following statement.

 \begin{assumption}
     The  $\zeta(t)$  is a Sub-Gaussian random vector. That is, $\exists \; \nu_1 > 0$ such that 
 $ \mathbb{E}\Big{[}  \exp{  \Big{(} \lambda
         \inprod{x}{\zeta(t)} \Big{)}}| \mathcal{F}_t \Big{]} \leq \exp{\Big{(} \frac{\norm{x}^2 \nu_1^2 \lambda^2}{2} \Big{)}} \; \; \forall \; x \in \mathbb{R}^n.$
     \label{ass4}
 \end{assumption}
Before proceeding, we present the following lemma, which is essential for proving the main theorem of this subsection. The proof of this lemma incorporates some steps from \cite{jin2019shortnoteconcentrationinequalities}.
\begin{lemma}
    There exists $\nu_2 \in \mathbb{R}^+$ and $a > 0$  such that 
    \vspace{-0.6em}
    \begin{equation*}
        \mathbb{E}\Big{[}  \exp{\Big{(} \lambda \norm{\zeta(t)}_\ast^2 \Big{)}} | \mathcal{F}_t \Big{]}  \leq \exp{( 0.5 \nu_2^2 \lambda^2 )} \; \; \forall \; |\lambda| \leq a^{-1}.
    \end{equation*}
    \label{lemmassub}
\end{lemma}
\vspace{-1em}
\begin{proof}
Consider the probability space $(\Omega, \mathcal{F} , \mathbb{P})$.  assume that for each  $\mathcal{F}_t$, there exists a regular conditional probability measure, $(\mathbb{Q}_{t,\omega})$ (see \cite{DU04} for more details) such that for any random variable $X$, we have
$ \mathbb{E}[X | \mathcal{F}] (\omega) = \int\limits_{\Omega} X d\mathbb{Q}_{t,\omega} = \mathbb{E}_{\mathbb{Q}_{t,\omega}} [X] (\omega)$.
 In this context, consider $\norm{\zeta(t)}_\ast = \sup\limits_{\norm{x} = 1} \inprod{x}{\zeta(t)}.$
 Since the constraint set is compact, we have $\norm{\zeta(t)}_\ast = \inprod{x}{\zeta(t)}$ where $\norm{x} = 1$.
Define the set $\mathcal{X} = \{x \in \mathbb{R}^n \; | \norm{x} = 1 \}$. Consider  an open cover of the set $\mathcal{X}$ as 
$\mathcal{X} \subseteq \cup_{x \in \mathcal{X}} U(x,\frac{1}{2})$. Since $\mathcal{X}$ is compact, there exists a finite subcover. In other words,
$\mathcal{X} \subseteq \cup_{i=1}^{C} U(x_i,\frac{1}{2})$. Hence,  $\exists \; x_i$ $i = 1, 2, \ldots, C$
\begin{equation*}
    \begin{split}
       & \norm{\zeta(t)}_\ast  =  \inprod{x}{\zeta(t)} = \inprod{x-x_i}{\zeta(t)}  + \inprod{x_i}{\zeta(t)}
        \\ & \leq \norm{x-x_i} \norm{\zeta(t)}_\ast  + \inprod{x_i}{\zeta(t)} \leq \frac{1}{2} \norm{\zeta(t)}_\ast + \inprod{x_i}{\zeta(t)}.   \end{split}
\end{equation*}
In other words,  for any $\delta > 0$, if  $\norm{\zeta(t)}_\ast \geq \delta$, then there exists some $  x_i$ (with $i \in [C])$ such that $\inprod{x_i}{\zeta(t)} \geq \frac{\delta}{2}$.  Since $\mathbb{Q}_{t, \omega}$ is a probability measure, we have
\\ $ \mathbb{Q}_{t,\omega}(\norm{\zeta(t)}_\ast \geq \delta) \leq \sum\limits_{i=1}^{C} \mathbb{Q}_{t,\omega}( \inprod{x_i}{\zeta(t)} \geq \frac{\delta}{2})$
\\ $= \sum\limits_{i=1}^{C}  \mathbb{Q}_{t,\omega} (\exp{  (\lambda \inprod{x_i}{\zeta(t)})} \geq \exp{(\lambda \frac{\delta}{2})})$
\\ $ \leq   \sum\limits_{i=1}^{C}\frac{\mathbb{E}_{\mathbb{Q}_{t,\omega}} [\exp{  (\lambda \inprod{x_i}{\zeta(t)})}]}{\exp{(\lambda \frac{\delta}{2})}} =  \sum\limits_{i=1}^{C} \frac{\mathbb{E}[\exp{  (\lambda \inprod{x_i}{\zeta(t)})}| \mathcal{F}_t]}{\exp{(\lambda \frac{\delta}{2})}} $
\\ $ \leq  C \exp \Big{(}{\Big{(}\frac{\lambda^2 \nu_1^2}{2} \Big{)} - \lambda \frac{\delta}{2}} \Big{)}$.
The first inequality of the above equation  follows  from the conditional Markov inequality. We can obtain a much tighter bound by taking the infimum on the RHS with respect to $\lambda$. Hence, we get  $ \mathbb{Q}_{t,\omega}(\norm{\zeta(t)}_\ast \geq \delta)  \leq C \exp{\Big{(} -\frac{\delta^2}{8 \nu_1^2} \Big{)}}$.

Hence, from \cite[Theorem 2.13]{Wainwright_2019}  we get that $\exists \; a > 0$
\begin{equation*}
\mathbb{E}_{\mathbb{Q}_{t,\omega}} \Big{[}\exp{\Big{(} \lambda \norm{\zeta(t)}_\ast^2 \Big{)}}\Big{]}  \leq \exp{\Big{(} \frac{\nu_2^2 \lambda^2}{2} \Big{)}} \; \; \; \text{such that} \; \forall \; |\lambda| \leq \frac{1}{a}
\end{equation*}
where, $\nu_2$ is a constant satisfying $\mathbb{E}_{\mathbb{Q}_{t,\omega}} \Big{[} \norm{\zeta(t)}_\ast^4 \Big{]} \leq \nu_2^2$. Note that this $\nu_2$ is independent of $\omega$. 
\end{proof}
Next, we provide the main results of this subsection. 
\begin{theorem}
    Consider any $\epsilon > 0$ then    $ \exists \; t_0 \in \mathbb{N}$ such that
    then under Assumption \ref{ass4} and $\forall \; t \geq t_0$, the following inequality holds
     $\mathbb{P}(f(z(t))-f^\ast \geq \epsilon) 
           \leq  \exp{\Bigg{(} \frac{-\epsilon^2 \mathrm{K}^2 (\sum\limits_{k=1}^{t} \alpha(k))^2}{18 D^2 \nu_1^2 \sum\limits_{k=1}^{t}\alpha(k)^2} \Bigg{)}}$
           \vspace{-1em}
    \begin{equation}
        \begin{split}
         +\frac{\exp{\Big{(} \sum\limits_{k=1}^{t} \frac{\alpha(k)^2}{2\sigma_R} \kappa_1 \Big{(} G^2 + B(t)^2 + \frac{\nu_2^2 \alpha(k)^2 \kappa_1}{4\sigma_R} \Big{)} \Big{)} }}{\exp \Big{(}{\frac{\epsilon}{3} \mathrm{K} \sum\limits_{k=1}^{t} \alpha(k)}\Big{)}} .
        \end{split}
        \label{gshs}
    \end{equation}
    \label{theorem5}
\end{theorem}
\vspace{-1.1em}
\begin{proof}
    To prove this theorem, we  consider $\eqref{24}$ in the proof of Theorem \ref{finite} and to bound the RHS of \eqref{24}, we utilize the Chernoff bound technique.
 Note that from Markov's inequality, $(\lambda > 0)$
 $  \mathbb{P} \Big(\sum\limits_{k=1}^{t} W(k) \geq \frac{\epsilon}{3} \mathrm{K} \sum\limits_{k=1}^{t} \alpha(k) ) \leq $ 
 \vspace{-0.9em}
 \begin{equation}
   \mathbb{P} (\exp ( \lambda \sum\limits_{k=1}^{t} W(k) ) \geq \exp ( \lambda  \tau(t)) ) \leq  \frac{\mathbb{E}\Big{[}\exp \Big{(}\lambda \sum\limits_{k=1}^{t} W(k) \Big{)} \Big{]}}{\exp (  \lambda \tau(t) )  }.
   \label{tete}
 \end{equation}
  To obtain a bound on the numerator of the RHS of \eqref{tete} , we observe that 
  \vspace{-0.7em}
  \begin{equation*}
      \begin{split}
          & \mathbb{E}\Big{[} \prod\limits_{k=1}^{t} \exp(\lambda W(k)) \Big{]} =  \mathbb{E}\Big{[} \mathbb{E}\Big{[} \prod\limits_{k=1}^{t} \exp( \lambda W(k)) | \mathcal{F}_t \Big{]}\Big{]}.
      \end{split}
  \end{equation*}
   From the definition it is clear  that $W(1)$, $W(2)$, $\ldots$ ,$W(t-1)$ are $\mathcal{F}_t$ measurable. Hence, we get
\vspace{-1em}
  \begin{equation*}
      \begin{split}
          &  \mathbb{E}\Big{[} \mathbb{E}\Big{[} \prod\limits_{k=1}^{t} \exp(\lambda W(k)) | \mathcal{F}_t \Big{]}\Big{]}
          \\ = & \mathbb{E} \Big{[}  \prod\limits_{k=1}^{t-1} \exp{(\lambda W(k))} \mathbb{E}[\alpha(t) a(t) \lambda \inprod{\zeta(t)}{x^\ast-x(t)} | \mathcal{F}_t] \Big{]}
      \end{split}
  \end{equation*}
  \vspace{-1.1em}
  \begin{flalign}
     \leq   \mathbb{E} \Big{[}  \prod\limits_{k=1}^{t-1} \exp{(\lambda W(k))} \Big{]}\exp{\Big{(}\frac{\alpha(t)^2 a(t)^2 \lambda^2 \nu_1^2 D^2}{2}\Big{)}}. &&
      \label{hdhdd}
  \end{flalign}
  The inequality in \eqref{hdhdd} follows from Assumption \ref{ass4}. Note that here, $a(t) = \prod\limits_{k=1}^{t} (1+ \rho(k))^{-1} < 1$. By continuing this process, we get  $\mathbb{E}\Big{[} \prod\limits_{k=1}^{t} \exp(W(k)) \Big{]} \leq \exp{(  2^{-1}\nu_1^2 D^2 \lambda^2  \sum\limits_{k=1}^{t}\alpha(k)^2 )}$.
 Hence, from \eqref{tete}, we get 
 \vspace{-0.75em}
 \begin{equation*}
 \mathbb{P} (\sum\limits_{k=1}^{t} W (k)\geq \tau(t))  \leq  \exp{( \lambda^2 \frac{\nu_1^2}{2} D^2 \sum\limits_{k=1}^{t}\alpha(k)^2 - \lambda \tau(t)} ).
 \end{equation*}
   By taking the infimum of the RHS of  w.r.t. $\lambda$, we get
   \vspace{-0.8em}
  \begin{equation}
      \begin{split}
          & \mathbb{P} (\sum\limits_{k=1}^{t} W(k) \geq \tau(t) ) \leq \exp {\Big{(}- \frac{\tau(t)^2}{2 D^2 \nu_1^2 \sum\limits_{k=1}^{t}\alpha(k)^2} \Big{)}}.
      \end{split}
      \label{rumia}
  \end{equation}
 Next, we need to obtain a bound on the right-hand side of the second term in \eqref{24}. Using the same procedure we get 
 \begin{flalign}
    \mathbb{P}(\sum\limits_{k=1}^{t} \frac{\alpha(k)^2}{2 \sigma_R} \norm{\Tilde{g}(k)}_\ast^2 \geq \tau(t) )       \leq  \frac{\mathbb{E} \Big{[} \exp{\Big{(} \sum\limits_{k=1}^{t} \frac{\alpha(k)^2}{2 \sigma_R} \norm{\Tilde{g}(k)}_\ast^2} \Big{)}\Big{]}}{\exp ({\tau(t)})}. && \nonumber
     \label{praneeth}
 \end{flalign}
Therefore, it is necessary to establish a bound for the numerator on the RHS. In doing so, we utilize Lemma \ref{lemmassub}. 
\\ $\mathbb{E} \Big{[} \exp{\Big{(} \sum\limits_{k=1}^{t} \frac{\alpha(k)^2}{2 \sigma_R} \norm{\Tilde{g}(k)}_\ast^2} \Big{)}\Big{]}$
\begin{equation*}
    \begin{split}
         = &  \mathbb{E} \Big{[} \prod\limits_{k=1}^{t-1} \exp{\Big{(}  \frac{\alpha(k)^2}{2 \sigma_R} \norm{\Tilde{g}(k)}_\ast^2} \Big{)} \mathbb{E} \Big{[} \exp{\Big{(}  \frac{\alpha(t)^2}{2 \sigma_R} \norm{\Tilde{g}(t)}_\ast^2 \Big{)}}| \mathcal{F}_t \Big{]}.
    \end{split}
\end{equation*}
Note that $\mathbb{E} \Big{[} \exp{\Big{(}  \frac{\alpha(t)^2}{2 \sigma_R} \norm{\Tilde{g}(t)}_\ast^2 \Big{)}}| \mathcal{F}_t \Big{]}$
\begin{flalign}
         \leq & \mathbb{E} \Big{[} \exp{\Big{(}  \frac{\alpha(t)^2}{2 \sigma_R} \kappa_1 \Big{(}\norm{g(t)}_\ast^2 + B(t)^2 +\norm{\Tilde{\zeta}(t)}_\ast^2 \Big{)} \Big{)}}| \mathcal{F}_t \Big{]} && \nonumber
        \\ \leq & \exp{\Big{(} \frac{\alpha(t)^2}{2 \sigma_R} \kappa_1 (G^2 + B(t)^2) \Big{)}} . \exp{\Big{(} \frac{\nu_2^2 \alpha(t)^4 \kappa_1^2 }{8 \sigma_R^2} \Big{)}}. &&
    \label{por}
\end{flalign}
The last inequality in \eqref{por} is due to Lemma \ref{lemmassub} where, with out loss of generality we have  assumed that $\frac{\alpha(t)^2 \kappa_1}{2\sigma_R} \leq \frac{1}{a}$.
By continuing this process we get  $\mathbb{E} \Big{[} \exp{\Big{(} \sum\limits_{k=1}^{t} \frac{\alpha(k)^2}{2 \sigma_R} \norm{\Tilde{g}(k)}_\ast^2} \Big{)}\Big{]}$
\\  $\leq  \exp{( \sum\limits_{k=1}^{t} \frac{\alpha(k)^2}{2\sigma_R} \kappa_1 \Big{(} G^2 + B(t)^2 + \frac{\nu_2^2 \alpha(k)^2 \kappa_1}{4\sigma_R} \Big{)} ) }$.
Hence, we get  $\mathbb{P}(\sum\limits_{k=1}^{t} \frac{\alpha(k)^2}{2 \sigma_R} \norm{\Tilde{g}(k)}_\ast^2 \geq \tau(t))$
\begin{flalign}
        \leq \exp{\Big{(} \sum\limits_{k=1}^{t} \frac{\alpha(k)^2}{2\sigma_R} \kappa_1 \Big{(} G^2 + B(t)^2 + \frac{\nu_2^2 \alpha(k)^2 \kappa_1}{4\sigma_R} \Big{)} - \tau(t)\Big{)} }. && \nonumber 
\end{flalign}
By adding the above equation with  \eqref{rumia},  \eqref{gshs} follows.
\end{proof}
 The next corollary  offers an alternative formulation of this theorem, illustrating how algorithmic convergence rates are affected by parameters like the bias $B(t)$. 
\begin{Corollary}
    For any $\epsilon > 0$ and a confidence level $0 < p < 1$, let $p_1 = 1-p$ and $\mathcal{R} = \sup\limits_{(x,y) \in \mathbb{X}} \mathbb{D}_R(x,y)$.

Define $t_0, t_1, t_2$  in such a way that $\forall \;  t \geq t_0$,

$\sum\limits_{k=1}^{t} \alpha(k) \geq \frac{3}{\epsilon} \exp{\Big{(}\sum\limits_{k=1}^{\infty} \frac{2 \alpha(k)B(k)}{\sigma_R}\Big{)}} \Big{(} \mathcal{R} + \sum\limits_{k=1}^{t} \alpha(k) B(k)\Big{)} ,$

similarly, $\forall \; t \geq t_1$ 
\\ $\Big{(}\sum\limits_{k=1}^{t} \alpha(k) \Big{)}^2 \geq \frac{18 D^2 \nu_1^2}{\epsilon^2} \ln \Big{(} \frac{2}{p_1} \Big{)} 
    \exp{\Big{(} \sum\limits_{j=1}^{\infty} \frac{4 \alpha(j) B(j)}{\sigma_R} \Big{)}}  \sum\limits_{k=1}^{t} \alpha(k)^2 $
\\ and $\forall \; t \geq t_2$
$ \sum\limits_{k=1}^{t} \alpha(k) \geq  \frac{3}{\epsilon}  \exp{\Big{(} \sum\limits_{j=1}^{\infty} \frac{2 \alpha(j) B(j)}{\sigma_R} \Big{)}} \times  \{ \ln{\Big{(}\frac{2}{p_1} \Big{)}}
   \\ +  \sum\limits_{k=1}^{t} \frac{\alpha(k)^2}{2 \sigma_R} \kappa_1 (G^2 + B(k)^2 + \frac{\nu_2^2 \kappa_1}{4 \sigma_R}) \}.$
Then $\forall \; t \geq t_0,t_1,t_2$ we have $f(z(t)) -f^\ast < \epsilon$ w.p. atleast $p$. 
\label{subgssh}
\end{Corollary}
Corollaries \ref{concentrationc} and \ref{subgssh} demonstrate that assuming a Sub-Gaussian stochastic gradient enhances the convergence rate, with terms $t_1$ and $t_2$ depending on $\ln{\frac{2}{1-p}}$ instead of $\frac{1}{1-p}$, while $t_0$ remains constant. Thus, selecting the right Bregman divergence is crucial.  This analysis highlights the significance of the bias term $B(t)$ in zeroth-order optimization, aiding in designing efficient algorithms and choosing optimal parameters.
\\ Before concluding, while the Sub-Gaussian assumption on noisy subgradients is more restrictive than standard assumptions, it holds in certain scenarios within zeroth-order optimization. For instance, if the function $f$ in \eqref{CP1} is continuously differentiable with a globally Lipschitz gradient, the approximated gradient from Gaussian approximation \cite{Nesterov2015RandomGM} satisfies Assumption \ref{ass4} (can be shown from \cite[Theorem-2.26]{Wainwright_2019}). We omit the detailed proof for brevity.
\vspace{-1.2em}
\section{Conclusion}
This letter provides insights into the SMD algorithm with non-zero diminishing bias, showing almost sure convergence to optimal solutions at a rate of $\mathcal{O}(\frac{1}{\sqrt{t}})$ under certain conditions on step-size.  Furthermore, we derived a finite-time analysis for the SMD algorithm with non-zero bias, elucidating how the algorithm's performance depends on various parameters, particularly step-size and bias. This analysis helps in understanding the practical aspects of algorithm performance over finite time periods. For future research, extending this analysis to other variants of the mirror descent algorithm, such as those with adaptive step-sizes or momentum terms, could further enhance its applicability.

\bibliographystyle{IEEEtran}
\bibliography{ref}


\end{document}